\documentclass[lettersize,journal]{IEEEtran}
\usepackage[caption=false,font=normalsize,labelfont=sf,textfont=sf]{subfig}
\usepackage{amsmath,amssymb}
\usepackage{natbib}
\setcitestyle{numbers,square}
\pdfoutput=1
\usepackage{color}
\usepackage{amsmath}
\usepackage{graphicx}
\usepackage{booktabs}
\usepackage{url}
\usepackage{multirow}
\usepackage{diagbox}
\usepackage{algorithm}
\usepackage{natbib}  
\newtheorem{remark}{Remark}
\newtheorem{theorem}{Theorem}
\newtheorem{assumption}{Assumption}
\newtheorem{lemma}{Lemma}
\newtheorem{proposition}{Proposition}
\newtheorem{corollary}{Corollary}
\newtheorem{example}{Example}
\newcommand{\beq}{\begin{equation}\begin{aligned}}
		\newcommand{\eeq}{\end{aligned}\end{equation}}
\newcommand{\beqn}{\begin{equation*}\begin{aligned}}
		\newcommand{\eeqn}{\end{aligned}\end{equation*}}
\begin{document}

\title{An Efficient Sparse Identification Algorithm For Stochastic Systems With General Observation Sequences}

\author{Ziming Wang, Xinghua Zhu*
	\thanks{Z. M. Wang is with NCMIS, LSEC, Academy of Mathematics and Systems Science, Beijing, P. R. China. Emails: wangziming@lsec.cc.ac.cn. X. H. Zhu is with the Key Laboratory of Systems and Control, Academy of Mathematics
		and Systems Science, Chinese Academy of Sciences, Beijing, P. R. China. Emails: zxh@amss.ac.cn.}
	\thanks{* Corresponding author.}
}

\markboth{IEEE TRANSACTIONS ON CYBERNETICS}%
{Shell \MakeLowercase{\textit{et al.}}: A Sample Article Using IEEEtran.cls for IEEE Journals}

\maketitle

\begin{abstract}                
This paper studies the sparse identification problem of unknown sparse parameter vectors in stochastic dynamic systems. Firstly, a novel sparse identification algorithm is proposed, which can generate sparse estimates based on least squares estimation by adaptively adjusting the threshold. 
Secondly, under a possibly weakest non-persistent excited condition, we prove that the proposed algorithm can correctly identify the zero and nonzero elements of the sparse parameter vector using a finite number of observations, and further estimates of the nonzero elements almost surely converge to the true values. Compared with the related works, e.g.,  LASSO, our method only requires the weakest assumptions and does not require solving additional optimization problems. Besides, our theoretical results do not require any statistical assumptions on the regression signals, including independence or stationarity, which makes our results promising for application to stochastic feedback systems.
Thirdly, the number of finite observations that guarantee the convergence of the zero-element set of unknown sparse parameters of the Hammerstein system is derived for the first time. 
 Finally, numerical simulations are provided, demonstrating the effectiveness of the proposed method. 
Since there is no additional optimization problem, i.e., no additional numerical error, the proposed algorithm performs much better than other related algorithms.
\end{abstract}


\section{Introduction}
Parameter estimation or filtering plays a very important role in system identification and control, signal processing, statistical learning, and other fields (\cite{blair2000multitarget,vaezi2015piecewise}). In recent decades, the classical identification and estimation theory has made great progress, and a relatively complete theoretical system has been formed (\cite{guo1990adaptive,chen1991identification}).

We note that in many practical application scenarios, such as the selection of effective basis functions in Hamiltonian systems, feature selection and robust inference of biomarkers based on omics data in personalized healthcare, channel estimation in ultra-wideband communication systems, since many components in the unknown parameter vector have no or negligible contribution to the system (these elements are zero or close to zero), the unknown parameter vector in a large number of systems is high-dimensional and sparse. Naturally, an interesting research direction is sparse parameter identification, i.e., the precise identification of zero and nonzero elements in the unknown parameter vector, in order to reduce model redundancy and obtain a leaner, better performance, and more reliable prediction model.

It is found that the research on sparse parameter identification has made considerable progress. In the field of signal processing, we find that the research of sparse parameter identification is based on compressed sensing (CS), i.e., solving an $L_0$(the number of nonzero elements) optimization problem with $L_2$(Euclidean metric) constraints. Since this optimization problem is NP-hard, researchers have developed some sparse signal estimation algorithms by transforming $L_0$ optimization into a more solvable $L_1$ convex optimization problem using the norm equivalence principle. For example, by using compressed sensing methods, Kalouptsidis et al. in \cite{KALOUPTSIDIS20111910} developed a sparse identification algorithms based on Kalman filtering and Expectation-Maximization, and verified the efficiency of the proposed algorithm via simulations. Cand{\`e}s and Tao et al. in \cite{candes2006stable,candes2006compressive,candes2005decoding} developed a relatively perfect signal reconstruction theory under some prior assumptions of the sparsity of the model. See \cite{kopsinis2010online,chen2009sparse,gu2009l_} for more references. The variable selection problem plays a very important role in the field of statistical learning. The well-known LASSO (the least absolute shrinkage and selection operator) is a linear regression method with $L_1$ regularization, and sparse identification algorithms represented by it and its variants have been extensively studied. Tibshirani in \cite{tibshirani1996regression} first proposed the LASSO in 1996. Zhao and Yu in \cite{zhao2006model} established the model selection consistency under regularity conditions. Zou in \cite{zou2006adaptive} first proposed the adaptive lasso and discussed the advantages of the adaptive lasso in detail. 
Combined with the ubiquitous sparsity in practical systems, we find that CS and adaptive lasso methodology are inherited and developed in the field of system identification and control, and the sparse parameter identification theory of stochastic dynamic systems is established. For example, Satheesh and Arun in \cite{PEREPU2015795} proposed a CS-based iterative basis pursuit denoising (BPDN) algorithm to estimate the unknown sparse parameters in the ARMAX model, where the number of zero elements of unknown parameters is known.
Roland et al. in \cite{6160383} studied the consistency of estimation for sparse linear regression models based on CS under the condition that the input signals are independent and identically distributed (i.i.d.).
Considering that the regressors in the stochastic dynamic system are often generated by the past input and output signals, so the independence assumption about the observation data is difficult to be satisfied. Based on this, Xie and Guo in \cite{xie2020analysis} developed a compressed consensus normalized least mean squares and provided the stability analysis of the algorithm based on CS theory, where the strong statistical assumptions such as independency are not required, but the prior assumption of the sparsity of the model is. Zhao et al. in \cite{zhao2020sparse} proposed a LASSO-based adaptive sparse estimation algorithm with general observation sequences. The zero elements in the unknown parameter are correctly identified with a finite number of observations and a non-persistent excited condition, which is the weakest excited condition in the existing literature about sparse parameter identification. There is also related research in some other fields, such as pruning in deep learning \cite{alzubi2020optimal,xu2020simultaneously}, etc.

Summarizing the existing research, we find that there are still some shortcomings in the sparse parameter identification theory. Firstly, the sparse estimation algorithm based on CS theory requires prior knowledge about the sparsity of the unknown parameter and the regression vectors. Secondly, the theoretical analysis of the LASSO-based sparse estimation algorithm requires a regularity condition of the regressor, which is a persistent excited condition actually and difficult or almost impossible to satisfy in stochastic feedback systems. Even the non-persistent excited condition for regressor in \cite{zhao2020sparse} is much stronger than the weakest excited condition proposed by Lai and Wei in \cite{lai1982least} (See subsection \ref{rr} and Table \ref{t111} for detail).
Thirdly, for a specific system, such as the Hammerstein system, it is very instructive to give a definite number of finite observations under further assumptions about the regressor, but this cannot be done in \cite{zhao2020sparse}.

Our contributions to the paper are summarized as follows: 
\begin{enumerate}
    \item We first propose a novel sparse parameter identification algorithm to estimate unknown and sparse parameters in stochastic regression systems. The principle of the algorithm is to generate a sparse estimate based on the least squares estimator by adjusting the threshold adaptively. Compared with the existing LASSO and its variants method for generating sparse estimates by optimizing the criterion function with a penalty term, our algorithm will not need to solve the optimization problem, so it will be more concise and efficient.
    \item Unlike the classical LS estimator's asymptotical theory (\cite{lai1982least}), under the same condition, we prove that the proposed algorithm can correctly identify the set of zero elements and nonzero elements in the unknown sparse parameter vector under finite observation, which is the convergence of the set. It is worth pointing out that this result is parallel to the convergence results in \cite{zhao2020sparse}, but the non-persistent excited condition required is much weaker than \cite{zhao2020sparse}. Furthermore, estimates of the nonzero will converge to the true values almost surely  with a convergence rate of $O\left(\sqrt{\frac{\log R_N}{\lambda^N_{\min}}}\right)$, i.e., parameter convergence, which, to the best of the authors' knowledge,  has never been achieved in sparse parameter identification algorithms for stochastic regression models.

\item For a classical system, such as the Hammerstein system, whose strong convergence and convergence rate are established by Zhao in \cite{5378462}, we give the number of finite observations under the same conditions, which is never achieved by sparse parameter identification in stochastic dynamic systems.

\item Simulation experiments compare the proposed algorithm with the classical least squares (LS) estimation, LASSO estimation, and algorithm in \cite{zhao2020sparse}. In contrast, other algorithms cannot exactly identify the element $0$, showing the great advantage of the algorithm proposed in this paper.
\end{enumerate}

The rest of this paper is organized as follows. In Section \ref{PROBLEM FORMULATION}, we give the problem formulation. 
Section \ref{Convergence} presents the main results of this paper, including the parameter convergence, the set convergence, and the comparison of conditions for consistency of other algorithm and the proposed algorithm in this paper. Section \ref{APPLICATION} applies the algorithm to the sparse parameter estimation of the Hammerstein system. A simulation example is given in Section \ref{simu}, and the concluding remarks are made in Section \ref{conclu}.

\section{Problem Formulation}\label{PROBLEM FORMULATION}
\subsection{Some preliminaries}\label{Preiminary}
 Let $(\Omega, \mathcal{F},\mathcal{P})$ be the probability space, $\omega$ be an element in $\Omega$, and $\mathbb{E}(\cdot)$ be the expectation operator. Denote $\|\cdot\|$ as the $2$-norm of vectors or matrices in this paper. For two positive sequences $\left\{a_k\right\}_{k \geq 1}$ and $\left\{b_k\right\}_{k \geq 1}$, $a_k=O\left(b_k\right)$ means $a_k \leq c b_k$ for $k \geq 1$ and some $c>0$, while $a_k=o\left(b_k\right)$ means $a_k / b_k \rightarrow 0$ as $k \rightarrow \infty$. 
\subsection{Sparse identification algorithm}\label{Distributed LS Algorithm}
Consider the parameter identification problem of the following discrete-time stochastic regression model,
\beq\label{equ1}
y_{k+1}=\varphi_k^T\theta +w_{k+1}, \quad k \geq 0
\eeq
where $y_{k+1}$ is the scalar observation or output at time $k$, $\varphi_k \in \mathbb{R}^r$ is the r-dimensional stochastic regression vector which may be the function of current and past inputs and outputs, $\theta\in \mathbb{R}^r$ is an unknown
$r$-dimensional parameter to be estimated, and $w_{k+1}$  is the stochastic noise sequence.

The above model \eqref{equ1} includes many parameterized systems, such as ARX system and Hammerstein system. We further denote the parameter vector $\theta$ and the index set of its zero elements by

\begin{align}
	\theta & \triangleq(\theta(1), \cdots, \theta(r))^T \notag\\
	\label{H}H^* & \triangleq\{l \in\{1, \cdots, r\} \mid \theta(l)=0\}.
\end{align}

The classical identification algorithms, such as the LS algorithm, can generate consistent estimates for the unknown parameters as the number of data goes to infinity. However, due to the existence of system noises and the finite number of observations in practice, it is hard to deduce from such estimates which parameters are exactly zero or not. In the existing literature, $L_1$ penalty term is added to the criterion function to be optimized to generate sparse estimates so the zero elements in the unknown parameter $\theta$ can be identified. However, there are two issues worth thinking about. Firstly, the existence of the penalty term destroys the optimality of LS and makes the estimation inevitably biased, so a stronger excited condition is needed to ensure the consistency of the algorithm. Secondly, solving the optimization problem with the penalty term brings great difficulty.
Our problem is to design a new sparse adaptive estimation algorithm that does not require penalty terms to infer the set $H^*$ in a finite number of steps and identify the unknown parameter $\theta$ by using stochastic regression vectors and the observation signals $\left\{\varphi_k, y_{k+1}\right\}_{k=1}^n$.
\begin{algorithm}[!ht]
	\caption{Sparse Identification Algorithm}\label{algorithm1}
	
	\textbf{Step 0:}(Initialization) Choose a positive sequence $\left\{\alpha_n\right\}$ satisfying $\alpha_n\to 0$ as $n\to \infty$ and
	$$
	\begin{aligned}
		\sqrt{\frac{\log R_n}{\lambda_{\min}^n}}=o\left(\alpha_n\right).
	\end{aligned}
	$$
	\textbf{Step 1:} Based on $\left\{\varphi_k, y_{k+1}\right\}_{k=1}^n$, begin with an initial vector $\theta_0$ and an initial matrix $P_0>0$, compute the matrix $P_{n+1}^{-1}$ and the estimate $\theta_{n+1}$ of $\theta$ by the following  recursive LS algorithm,\\
	\begin{align}\label{e:sf1}
		P_{k+1}&=P_{k}-a_{k}P_{k}\varphi_k\varphi^T_kP_k,\\
		\label{e:sf12}a_k&=\frac{1}{1+\varphi^T_kP_k\varphi_k},\\
		\label{e:sf13}\theta_{k+1}&=\theta_k+a_kP_k\varphi_k\big\{y_{k+1}-\varphi^T_k\theta_k\big\}.
	\end{align}
	\textbf{Step 2:} Obtain $\beta_{n+1}(l)$ by the following mechanism,
	\begin{align}\label{e:sf2}
		\beta_{n+1}(l)=\left\{
		\begin{array}{ll}
			0\quad if\ |\theta_{n+1}(l)|<\alpha_{n+1}\\
			\theta_{n+1}(l) \quad else
		\end{array}
		\right.,
	\end{align} 
	then obtain 
	\begin{align}\label{e:sf3}
		\beta_{n+1}&=\left(\beta_{n+1}(1), \cdots, \beta_{n+1}(r)\right)^T \\
		\label{Hn}H_{n+1} & \triangleq\left\{l=1, \cdots, r \mid \beta_{n+1}(l)=0\right\}.
	\end{align}
	
\end{algorithm}

\section{Theoretical Properties of the Sparse Identification Algorithm}\label{Convergence}
In this section, we will investigate the asymptotic analysis of the unknown sparse parameter vector and the convergence of the sets of zero elements  with a finite number of observations, 
To this end, we introduce the following assumptions to be used for the theoretical analysis.
\begin{assumption}\label{c1} The noise $\left\{w_k, \mathcal{F}_k\right\}_{k \geq 1}$ is a martingale difference sequence, i.e., $\mathbb{E}\left[w_{k+1} \mid \mathcal{F}_k\right]=0, k \geq 1$, and there exists some $\gamma>2$ such that $\sup _k \mathbb{E}\left[\left|w_{k+1}\right|^\gamma \mid \mathcal{F}_k\right]<\infty$ a.s.
\end{assumption}
\begin{assumption} [Non-Excited Condition] \label{c3} 
	The growth rate of $\ \log \big(\lambda_{\max}\big\{P_{n+1}^{-1}\big\}\big)$ is slower than $\lambda_{\min}\big\{P_{n+1}^{-1}\big\}$, that is
	$$\lim\limits_{n\to \infty} \frac{\log R_n}{\lambda^n_{\min}} =0\ \ \ \ a.s.$$ 
	holds, where $R_n=1+\sum_{k=0}^n\|\varphi_k\|^2$,  $\label{xcv}
	\lambda_{\min}^n=\lambda_{\min}\{ P_0^{-1}+\sum_{k=0}^n  \varphi_k\varphi_k^T\}.$
\end{assumption}
\begin{remark}
A further explanation of this non-excited condition is provided in subsection \ref{rr}. 
\end{remark}

\subsection {Set and parameter convergence of estimates}
Assume that there are $d$ ($d$ is unknown) nonzero elements in vector $\theta$ . Without loss of generality, we assume $\theta=[\theta(1), \ldots, \theta(d), \theta(d+1), \ldots \theta(r)]^T$ and $\theta(i) \neq 0, i=1, \ldots, d, \theta(j)=0, j=d+1, \ldots, r$.
Before giving the main results, we first state a classical result in stochastic adaptive control.

\begin{lemma} [\cite{GUO1995435}]\label{lemma0-lasso}
	Assume that Assumptions \ref{c1} holds. Then as $n \rightarrow \infty$, the estimation error of the recursive LS algorithm \eqref{e:sf1}-\eqref{e:sf13} is bounded by
	$$
	\left\|\theta_{n+1}-\theta\right\|^2\leq C_0\frac{\log R_n}{\lambda^n_{\min}} \text { a.s., }
	$$
	where $C_0$ is a finite constant.
\end{lemma}
	\begin{remark}
The specific value of $C_0$ can be referred to \cite{GUO1995435}.
\end{remark}
\begin{remark}
	We remark that the recursive LS algorithm has strong consistency  under  Assumption \ref{c3}, i.e.,  $|\theta_{n+1}(l)-\theta(l)|\to 0\ a.s.$ for all $l=1,\cdots,r.$
\end{remark}
For the estimate $\beta_{n+1}$ generated by algorithms \eqref{e:sf1}-\eqref{e:sf3}, we have the following main results.
\begin{theorem}[Parameter Convergence]\label{lemma1-lasso}
	 Assume that Assumptions \ref{c1} - \ref{c3} hold. Then
	$$
	\beta_{n+1}(l) \underset{n \rightarrow \infty}{\longrightarrow} \theta(l), \quad l=1, \ldots, r \text { a.s. }
	$$
\end{theorem}
\begin{IEEEproof}
	By the definition of $\beta_{n+1}$, we have for $\forall l \in \{1,\ldots,r\}$
	\beqn
	|\beta_{n+1}(l)-\theta(l)| & \leq|\beta_{n+1}(l)-\theta_{n+1}(l)|+|\theta_{n+1}(l)-\theta(l)| \\
	&\leq \alpha_{n+1}+|\theta_{n+1}(l)-\theta(l)|.
	\eeqn
	Thus the proof is completed by the fact that $\alpha_n \to 0$ as $n \to \infty$ and the convergence results $|\theta_{n+1}(l)-\theta(l)|\to 0\ \ \ a.s.$
\end{IEEEproof}
\begin{theorem}[Set Convergence]\label{Set Convergence}
	 Assume that Assumptions \ref{c1} - \ref{c3} hold. Then there exists an $\omega$-space $\Omega_0$ with $\mathcal{P}\left\{\Omega_0\right\}=1$ such that for any $\omega \in \Omega_0$, there exists an integer $N_0(\omega)$ such that
	$$
	\beta_{N+1}(d+1)=\cdots=\beta_{N+1}(r)=0, \quad N \geq N_0(\omega),
	$$
	i.e., $H_{N+1}=H^*$, where $H_{N+1}$ and $H^*$ are defined by \eqref{Hn} and \eqref{H} respectively.
\end{theorem}
\begin{IEEEproof}
	Noting that Assumption \ref{c3} holds almost surely, there exists an $\Omega_0$
	with $\mathcal{P}\{\Omega_0\} = 1$ such that Assumption \ref{c3} holds for any $\omega \in \Omega_0$. In the
	following, we will consider the estimate sequence $\{\beta_{n+1}\}$ on a
	fixed sample path $\omega \in \Omega_0$.
	
	Firstly, we have for $i \in \{1,\ldots,d\}$,
	
	\beqn
	|\beta_{n+1}(i)|\geq&|\theta(i)|-|\beta_{n+1}(i)-\theta(i)|\notag\\
	\geq&|\theta(i)|-|\beta_{n+1}(i)-\theta_{n+1}(i)|-|\theta_{n+1}(i)-\theta(i)|\\
	\geq& |\theta(i)|-\alpha_{n+1}-|\theta_{n+1}(i)-\theta(i)|
	.
	\eeqn
	By the fact that $|\theta(i)|>0$ for $i\in \{1,
	\ldots,d\}$, $\alpha_n \to 0$ as $n \to \infty$ and the convergence results $|\theta_{n+1}(i)-\theta(i)|\to 0$,
	there exists a constant $N_1$ such that $\forall N>N_1$, 
	\beqn
	\alpha_{N+1}+|\theta_{N+1}(i)-\theta(i)|\leq\frac{|\theta(i)|}{2},
	\eeqn
	thus, we have for $\forall i \in \{1,\ldots,d\}$
	\beqn
	|\beta_{N+1}(i)|\geq\frac{|\theta(i)|}{2}>0.
	\eeqn
	
	Then for $i \in \{d+1,\ldots,r\}$, since $\theta(i)=0$, by Lemma \ref{lemma0-lasso} we have
	\beqn
	|\theta_{n+1}(i)|= O\left(\sqrt{\frac{\log R_n}{\lambda^n_{\min}}}\right) \text { a.s. }
	\eeqn
	By using the property $\sqrt{\frac{\log R_n}{\lambda_{\min}^n}}=o\left(\alpha_n\right)$, there exists $N_2$ such that, for $\forall N>N_2$ we have
	$$
	|\theta_{N+1}(i)|<\alpha_{N+1}.
	$$
	Then by the definition of $\beta_{N+1}$ we have $\beta_{N+1}(i)=0$ for $\forall N>N_2$ and $i\in \{d+1,\ldots,r\}$.
	
	Finally, the proof is completed by setting $N_0(\omega)=\max(N_1,N_2)$.
\end{IEEEproof}
\begin{remark}
Theorem \ref{Set Convergence} shows that the index set of the zero elements in $\theta$ can be correctly identified with a finite number of observations, and estimates for the nonzero elements will also be nonzero in a finite number of observations.
\end{remark}
The following corollary follows immediately from Theorem \ref{Set Convergence} and Lemma \ref{lemma0-lasso}.
\begin{corollary}\label{con1}
Under assumptions of Theorem \ref{Set Convergence}, Algorithm 1 has the following convergence rate as $N\to\infty$:
\beqn
\|\beta_{N+1}-\theta\|=O\left(\sqrt{\frac{\log R_N}{\lambda^N_{\min}}}\right)\quad a.s.
\eeqn
\end{corollary}
To the best of the authors' knowledge, this convergence rate is parallel to that of \cite{lai1982least,GUO1995435} and has never been achieved in sparse parameter identification algorithms for stochastic regression models.
Next, we compare the conditions of observation data that guarantee the convergence of Algorithm 1, LASSO, and its variants respectively.
\subsection{Comparison of conditions for consistency of LASSO, adaptive LASSO, the algorithm proposed in \cite{zhao2020sparse} and Assumption \ref{c3}}\label{rr}

\begin{table*}[htbp]\label{t111}
    \centering
	\caption{Conditions for consistency of LASSO and its variations and Algorithm \ref{algorithm1}}
	\label{tab:sparse_excitation}
	\tabcolsep=0.5cm
    \begin{tabular}{ll}
			\toprule
			 & Conditions on system\\
			\midrule
			LASSO \cite{zhao2006model} & Regularity condition: $D_n\to D$ as $n\to \infty$ \\
			& Strong irrepresentable condition: for some $\eta>0$, $|D_n^{21}(D_n^{11})^{-1}sgn(\theta_1)|\leq 1-\eta$ \\
			\midrule
			Adaptive LASSO \cite{zou2006adaptive} & Regularity condition: $D_n\to D$ as $n\to \infty$ \\
			\midrule
			\cite{zhao2020sparse} & $\frac{R_n}{\lambda_{min}^n}\sqrt{\frac{\log R_n}{\lambda_{min}^n}} \to 0$ as $n\to \infty$\\
			\midrule
			Algorithm \ref{algorithm1}& $\frac{\log R_n}{\lambda_{min}^n} \to 0$ as $n\to \infty$\\
			\bottomrule
	\end{tabular}
\end{table*}

For simplicity of notations, we still assume that the parameter vector $\theta=\left[\begin{array}{ll}\theta_1^T & \theta_2^T\end{array}\right]^T, \theta_1=[\theta(1) \ldots \theta(d)]^T, \theta_2=[\theta(d+1) \ldots \theta(r)]^T$ such that $\theta(i) \neq 0, i=1, \ldots, d$ and $\theta(j)=0, j=d+1, \ldots, r$.
Denote
$$
D_n \triangleq \frac{1}{n} \sum_{k=1}^n \varphi_k \varphi_k^T=\left[\begin{array}{ll}
D_n^{11} & D_n^{12} \\
D_n^{21} & D_n^{22}
\end{array}\right]
$$
where $D_n^{11} \in \mathbb{R}^{d \times d}$ and $D_n^{12}, D_n^{21}$ and $D_n^{22}$ are with compatible dimensions. The comparison on conditions for consistency of LASSO and its variations as well as Algorithm \ref{algorithm1} is made in Table \ref{tab:sparse_excitation}.
The strong irrepresentable condition given in Table \ref{tab:sparse_excitation} 
is actually prior structural information on the sparsity of the parameter vector,  which is not required in our paper. Furthermore, from Table \ref{tab:sparse_excitation}, we can directly verify that Assumption \ref{c3} includes the regularity condition as its special case.
Besides, The algorithm proposed by \cite{zhao2020sparse} requires an excited condition that $R_n\sqrt{\log R_n}=o\left((\lambda_{min}^n)^{3/2}\right)$, which fails for $R_n=O\left((\lambda_{min}^n)^{\delta}\right)$, $\forall \delta>3/2$.
However, from Theorem \ref{Set Convergence} and Corollary \ref{con1} we know that Algorithm \ref{algorithm1} can achieve the consistency for  $R_n=O\left((\lambda_{min}^n)^{\delta}\right)$, $\forall\delta>0$, which greatly improves the applicability of the algorithm.
As far as the authors know, Assumption \ref{c3} is the weakest excited condition compared with the one required for sparse parameter identification in the existing literature.

\section{Practical applications}\label{APPLICATION}
Many applications of the sparse identification algorithm have been shown in \cite{zhao2020sparse}, such as the identification of Hammerstein systems \cite{eskinat1991use} and linear stochastic systems with self-tuning regulation control \cite{GUO1995435}. However, \cite{zhao2020sparse} only gave the existence of the constant $N_0(\omega)$ as in Theorem \ref{Set Convergence}, but did not give its specific range, which is not conducive to the development of practical applications.
In this section, we take the identification of Hammerstein systems as an example, and provide the specific range of the constant $N_0(\omega)$.

Hammerstein system is a modular nonlinear dynamic system with a static nonlinear function followed by a linear dynamic subsystem. Due to its simple structure, as well as a good simulation of nonlinear, Hammerstein system has been widely used in many practical engineering scenarios, such as power amplifiers, manipulators, as well as the neutralization reaction in chemical processes, lithium-ion battery heating systems and heat exchange system, (cf., \cite{SMITH2007551},\cite{Eskinat1991255}), etc.

Considering a Hammerstein system whose linear subsystem is an ARX system and whose nonlinear function is a combination of basis functions:
\beq\label{e:h1}
y_{k+1}=& a_1 zy_{k+1}+\cdots+a_p z^py_{k+1} \\
&+b_1 f\left(u_k\right)+\cdots+b_q f\left(z^{q-1}u_k\right)+w_{k+1},
\eeq
\beq\label{e:h2}
f\left(u_k\right)= \sum_{j=1}^m c_j g_j\left(u_k\right)
\eeq
where $\left\{g_j(\cdot)\right\}_{j=1}^m$ are the basis functions, $z$ is the backward-shift operator (i.e., $zy_{k+1}=y_k$), $y_{k+1}$ and $u_k$ are the output and input, $w_{k+1}$ is the system noise, $a_1,\cdots,a_p,b_1,\cdots,$ $b_q,c_1,\cdots,c_m$ are unknown coefficients that need to be estimated. 
Denote
$$
\begin{array}{c}
\theta\triangleq\left(a_1 \ldots a_p\left(b_1 c_1\right) \ldots\left(b_1 c_m\right) \ldots\left(b_q c_1\right) \ldots\left(b_q c_m\right)\right)^T \\
\varphi_k\triangleq\left(y_k \ldots y_{k+1-p} g_1\left(u_k\right) \ldots g_m\left(u_k\right) \ldots\right. \\
\left.\quad g_1\left(u_{k+1-q}\right) \ldots g_m\left(u_{k+1-q}\right)\right)^T,
\end{array}
$$
the Hammerstein system can be written in the form of \eqref{equ1}. 

   In order to obtain a good approximation of nonlinear functions $f(\cdot)$, a large number of basis functions (which leads to a  very large $m$) are often used in practice, which will lead to the unknown parameter vector $\theta$ is probably to be high dimensional and sparse in practice. To obtain a simpler but more precise model of the system, we find that it is crucial to determine the number of the effective basis functions in $\left\{g_j(\cdot)\right\}_{j=1}^m$, or sparse identification of the unknown parameter vector $\theta$. Thus Algorithm 1 can be applied to infer the zero elements in $\theta$. 
   Denote
$$M=\left[\begin{array}{lll}
M(1) & \cdots & M(m)
\end{array}\right]
$$
with $M(l)=\left[b_1 c_l \cdots b_q c_l\right]^T,$ $l=1, \ldots, m$. Thus we can find that the noneffective basis functions in $\left\{g_j(\cdot)\right\}_{j=1}^m$ correspond to zero columns in matrix $M$. 
   Consequently,  executing algorithms \eqref{e:sf1}-\eqref{e:sf13} by using  the inputs and  outputs $\{u_k, y_{k+1}\}_{k=1}^n$ of Hammerstein system \eqref{e:h1}, we can obtain 
    $\theta_{n+1}$, followed by executing \eqref{e:sf2}-\eqref{Hn}, we have 
   \begin{align}
   &\beta_{n+1}=\left[a_{1, n+1} \ldots a_{p, n+1}\left(b_1 c_1\right)_{n+1} \ldots\left(b_1 c_m\right)_{n+1}\right. \notag\\ 
   &\qquad\left.\ldots\left(b_q c_1\right)_{n+1} \ldots\left(b_q c_m\right)_{n+1}\right]^T, \\ \label{M}
   &M_{n+1}={\left[M_{n+1}(1) \cdots M_{n+1}(m)\right], } \\ 
   &\text{with}\ M_{n+1}(l)=\left[\left(b_1 c_l\right)_{n+1} \ldots\left(b_q c_l\right)_{n+1}\right]^T, l=1, \ldots, m,\notag \\ 
   &H^*=\left\{l\in\{1, \ldots, m\} \mid c_l=0\right\}, \\ 
   &H_{n+1}=\left\{l=1, \ldots, m \mid M_{n+1}(l)=0\right\} . 
   \end{align}
Before presenting the main results, we need the following assumptions. 
\begin{assumption}\label{b1}
$\left\{1, g_1(x), \ldots, g_m(x)\right\}$ is linearly independent over some interval $[a, b]$.
\end{assumption} 
\begin{assumption}\label{b2} 
Polynomial $A(z)=1-a_1 z-\cdots-a_p z^p$ is stable, i.e., $|A(z)| \neq 0, \forall|z| \leq 1$ and $b_1^2+\cdots+b_q^2 > 0$.
\end{assumption} 
\begin{assumption}\label{b3} 
$\left\{u_k\right\}_{k \geq 1}$ is an i.i.d. sequence with density $p(x)$ which is positive and continuous on $[a, b]$ and $0<\mathbb{E} g_j^2\left(u_k\right)<\infty, j=$ $1, \ldots, m$. Further, $\left\{u_k\right\}_{k \geq 1}$ and $\left\{w_k\right\}_{k \geq 1}$ are mutually independent.
\end{assumption} 
\begin{remark}
Assumptions \ref{b1}-\ref{b3} are essential to ensure convergence of parameter estimates in the Hammerstein system, as described in detail in \cite{5378462}.
\end{remark}
\begin{lemma} [\cite{5378462}] If Assumptions \ref{c1} and \ref{b1}-\ref{b3} hold, then for the maximal and minimal eigenvalues of $\sum_{k=1}^n \varphi_k \varphi_k^T$, the following inequalities
\begin{align}\label{c1c2c3c4}
C_1 n \leq& R_n \leq C_2 n, \text { a.s. } \\\label{c1c2c3c41}
C_3 n \leq& \lambda^n_{\min} \leq C_4 n \text {, a.s. }
\end{align}
hold for some $0<C_1<C_2, 0<C_3<C_4$. Further, the LS estimate $\theta_{n+1}$ has the following convergence speed,
$$
\left\|\theta_{n+1}-\theta\right\|=O\left(\sqrt{\frac{\log n}{n}}\right) \text { a.s. }
$$
\end{lemma}

\begin{proposition} 
Set $\alpha_n=\{\frac{\log n}{n}\}^{\epsilon}$ for any fixed $\epsilon \in\left(0, \frac{1}{2}\right)$. If Assumptions \ref{c1} and \ref{b1}-\ref{b3} hold for Hammerstein system \eqref{e:h1}-\eqref{e:h2}, then there exists an $\omega$-set $\Omega_0$ with $\mathcal{P}\left\{\Omega_0\right\}=1$ such that for any $\omega \in \Omega_0$,
$$
H_{N+1}=H^* \quad \forall N \geq N_0(\omega)
$$
holds for $N_0(\omega)$ being an positive integer,
i.e., the effective basis functions in $\left\{g_j(\cdot)\right\}_{j=1}^m$ can be correctly identified.
\end{proposition}
\begin{IEEEproof} By \eqref{c1c2c3c4}-\eqref{c1c2c3c41} and noticing $\alpha_n=\{\frac{\log n}{n}\}^{\epsilon}, \epsilon \in\left(0, \frac{1}{2}\right)$, we can verify that Assumptions \ref{c1}-\ref{c3} hold for the regression model composed of \eqref{e:h1}-\eqref{e:h2}, thus, the desired results followed by Theorem \ref{Set Convergence} directly.
\end{IEEEproof}

To give a specific selection method of the positive integer $N_0(\omega)$, we need the following assumption.
\begin{assumption}\label{c4}
$|\theta(i)|>C_5>0$ for $i=1,\ldots,d$.
\end{assumption}
\begin{remark}This assumption is necessary and natural, and if $|\theta(i)|>C_5$ $i=1,\ldots,d$ is very close to 0, the required constant $N_0(\omega)$ is inevitably going to be very large.
\end{remark}

Without loss of generality, set the threshold value $\alpha_n$ in the form of  $\alpha_n=M\left(\frac{R_n}{\lambda_{min}^{n}}\right)^{\epsilon}$ with $M>0$ and $0<\epsilon<\frac{1}{2}$. 
\begin{theorem}\label{sparsethm:specific}
Under Assumptions \ref{c1}-\ref{c4}, there exists an $\omega$-set $\Omega_0$ with $\mathcal{P}\left\{\Omega_0\right\}=1$ such that for any $\omega \in \Omega_0$,
\beq\label{eq:maxn}
H_{N+1}=H^*\quad \forall N>N_0
\eeq holds, where 
\begin{align} N_0=&\max\left\{\frac{47}{C_2},\frac{2k_1}{C_3}\log\left(\frac{C_2k_1}{C_3}\right),\frac{2k_2}{C_3}\log\left(\frac{C_2k_2}{C_3}\right)\right\},\notag\\\label{N_0}
k_1=&\left(\frac{\sqrt{C_0}}{M}\right)^{\frac{2}{1-2\epsilon}}, k_2=\left(\frac{2M}{C_5}\right)^{\frac{1}{\epsilon}}. 
\end{align}
\end{theorem}
\begin{IEEEproof}
From the proof of Theorem \ref{Set Convergence} we know that, the following two formulas are sufficient conditions for $H_{N+1}=H^*$:
\begin{align}
|\theta_{N+1}(i)-\theta(i)|+\alpha_{N+1}<|\theta(i)|, \quad i\in\{1,\ldots,d\}, \label{eqhkzhkk}\\
|\theta_{N+1}(i)|<\alpha_{N+1},\quad i\in\{d+1,\ldots,r\}. \label{eqhkzhkz}
\end{align}
Thus, if we can verify \eqref{eqhkzhkk} and \eqref{eqhkzhkz} under the above $N$, then the desired result \eqref{eq:maxn} is true naturally.
We first provide an inequality as follows,
$$
\frac{\log N}{N}<\frac{1}{t},\quad \forall N>\max\{47,2t\log t\},\ t>0.
$$
In fact, for $0<t<10$, we have
$$
\frac{\log N}{N}<\frac{\log 47}{47}<\frac{1}{10}<\frac{1}{t}.
$$
For $t\geq 10$, we have
$$
\frac{\log N}{N}<\frac{\log t + \log 2 +\log\log t}{2t\log t}<\frac{2\log t}{2t\log t}=\frac{1}{t}.
$$
On one hand, by \eqref{eq:maxn} and \eqref{N_0} we have
$$
C_2N>\max\left\{ 47, \frac{2C_2k_1}{C_3}\log\left(\frac{C_2k_1}{C_3}\right) \right\},
$$
thus $\frac{\log(C_2 N)}{C_2 N}<\frac{C_3}{C_2}\left(\frac{M}{\sqrt{C_0}}\right)^{\frac{2}{1-2\epsilon}}$, that is $\frac{\log(C_2 N)}{C_3 N}<\left(\frac{M}{\sqrt{C_0}}\right)^{\frac{2}{1-2\epsilon}}$.

Now by Lemma \ref{lemma0-lasso} and \eqref{c1c2c3c4} as well as \eqref{c1c2c3c41} we have 
\beq\label{eq:hkkhkz}
\frac{|\theta_{N+1}(i)-\theta(i)|}{\alpha_{N+1}}&\leq\frac{\sqrt{C_0}}{M}\left(\frac{\log (C_2 N)}{C_3 N}\right)^{\frac{1}{2}-\epsilon}\\
&< \frac{\sqrt{C_0}}{M}\left( \frac{M}{\sqrt{C_0}}\right)^{\frac{1-2\epsilon}{1-2\epsilon}} \\
&=1.
\eeq
Thus \eqref{eqhkzhkz} is satisfied by the fact that $\theta(i)=0$ for $i=d+1,\ldots,r$. On the other hand, by 
$$
C_2N>\max\{ 47, \frac{2C_2k_2}{C_3}\log\left(\frac{C_2k_2}{C_3}\right) \},
$$           
we have $\frac{\log(C_2 N)}{C_3 N}<\left(\frac{C_5}{2M}\right)^{\frac{1}{\epsilon}}$.
Then by \eqref{c1c2c3c4} we have,
\beq\label{eq:hkkhkk}
2\alpha_{N+1}&\leq 2M\left(\frac{\log(C_2 N)}{C_3 N}\right)^{\epsilon} \\
&< 2M \left(\frac{C_5}{2M}\right)^{\frac{\epsilon}{\epsilon}} \\
&=C_5.
\eeq
\eqref{eqhkzhkk} is satisfied by combining \eqref{eq:hkkhkz} and \eqref{eq:hkkhkk}. This completes the proof of the theorem.
\end{IEEEproof}

In Theorem \ref{sparsethm:specific}, the constant $N_0$ is corresponding to $k_1$ and $k_2$. To quickly identify the position of the zero elements, we need to make $\max\{k_1,k_2\}$ as small as possible.
Clearly, for fixed $\epsilon$, $k_1$ decreases as $M$ increases, while $k_2$ increases as $M$ increases.
A straightforward calculation shows that
\beq\label{eqzkhzkh}
&\max\{k_1,k_2\} \geq k_1 \geq \frac{4C_0}{C_5^2},\quad \text{if } M\leq 2^{2\epsilon-1}C_0^{\epsilon}C_5^{1-2\epsilon},\\
&\max\{k_1,k_2\} \geq k_2 \geq \frac{4C_0}{C_5^2},\quad \text{if } M\geq 2^{2\epsilon-1}C_0^{\epsilon}C_5^{1-2\epsilon}.
\eeq
Thus the lower bound of $\max\{k_1,k_2\}$ is $4C_0C_5^{-2}$, which motivates the following corollary.
\begin{corollary}\label{sparsecol:specific}
Under conditions of Theorem \ref{sparsethm:specific}, with $M=2^{2\epsilon-1}C_0^{\epsilon}C_5^{1-2\epsilon}$, we have
\beq\label{eq:maxn2}
H_{N+1}=H^*,\ \forall N>\max\left\{\frac{47}{C_2},\frac{8C_0}{C_3C_5^2}\log\left(\frac{4C_2C_0}{C_3C_5^2}\right)\right\}.
\eeq
\end{corollary}

\begin{IEEEproof}
The proof is a combination of \eqref{eqzkhzkh} and Theorem \ref{sparsethm:specific}.
\end{IEEEproof}
\begin{remark} Noticing that
$
M \triangleq\left[b_1 \cdots b_q\right]^T\left[\begin{array}{lll}
c_1 & \cdots & c_m
\end{array}\right] \text {, }
$
we can further infer the estimates for the nonzero elements in $\left\{b_i, i=1, \ldots, q\right\}$ and $\left\{d_l, l=1, \ldots, m\right\}$ by using a singular value decomposition (SVD) algorithm to $M_{N+1}$ defined by \eqref{M}, see Chaoui et al. (2005) for more details.
\end{remark}
\section{Simulation Results }\label{simu}
In this section, we provide an example to illustrate the performance of the novel sparse identification algorithm (i.e., Algorithm 1) proposed in this paper.

\begin{example}\label{e1}
Consider a discrete-time stochastic regression model \eqref{equ1} with the dimension $r=10$. The noise sequence $\{w_{k+1}\}, k \geq 0$ in \eqref{equ1} is independent and identically distributed with $w_{k+1} \sim \mathcal{N}(0,0.1)$ (Gaussian distribution with zero mean and variance 0.1). Let the regression vector $\varphi_k \in$ $\mathbb{R}^r(k \geq 1)$ be generated by the following state space model,
\beqn
x_k &= A x_{k-1} + \varepsilon_k, \\
\varphi_k &= B_k x_k,
\eeqn
where $x_k \in \mathbb{R}^r$ is the state of the above system with $x_{0}=[\underbrace{1, \ldots, 1}_r]^T$, the matrices $A, B_k$ and vector $\varepsilon_k$ are chosen according to the following way,
$$
\begin{array}{l}
A_k=\operatorname{diag}\{\underbrace{1.01, \cdots, 1.01}_r\}, \\
B_k=\left\{b_{k,ij}\right\} \in \mathbb{R}^{r\times r} \text{ with } b_{k,ij}\sim \mathcal{N}(0,1), \\
\varepsilon_k=\left\{\varepsilon_{k,i}\right\} \in \mathbb{R}^{r} \text{ with } \varepsilon_{k,i}\sim \mathcal{N}(0,1).
\end{array}
$$

The true value of the unknown parameter is 
\beqn
\theta &= \left[\theta(1),\theta(2),\theta(3),\theta(4),\theta(5),\theta(6),\theta(7),\theta(8),\theta(9),\theta(10)\right]^T \\
&= \left[0.8,1.6,-0.3,0.05,0,0,0,0,0,0\right]^T.
\eeqn
Ten simulations are performed. Fig. \ref{fig1:sparse} shows the estimate sequences $\left\{\theta_n(1),\theta_n(2),\ldots,\theta_n(10)\right\}_{n=1}^{600}$ generated by Algorithm \ref{algorithm1} with $\alpha_n=0.1\left(\frac{R_n}{\lambda_{min}^{n}}\right)^{\frac{1}{4}}$ from one of the simulations. Table \ref{tab1:sparse} compares the averaged estimates from the ten simulations generated by least squares, LASSO with $\lambda_n=$ $n^{0.75}$ (see \cite{zhao2006model}), the algorithm proposed by \cite{zhao2020sparse} with $\lambda_n= (\lambda_{min}^n)^{0.75}$ and Algorithm \ref{algorithm1} for $\theta(5),\theta(6),\theta(7),\theta(8),\theta(9),\theta(10)$, with different data length $n$. 
We adopt the Python CVX tools (http://cvxopt.org/) to solve the convex  optimization in LASSO (\cite{zhao2006model}) and the algorithm proposed by \cite{zhao2020sparse}.

From Fig. \ref{fig1:sparse} and Table \ref{tab1:sparse}, we can find that,
compared with the least squares estimates, the LASSO (\cite{zhao2006model}), the algorithm proposed by \cite{zhao2020sparse} and Algorithm \ref{algorithm1} generate sparser and more accurate estimates for the system parameters, naturally, give us valuable information in inferring the zero and nonzero elements in the unknown parameters. 
Moreover, compared with LASSO (\cite{zhao2006model}) and the algorithm proposed by \cite{zhao2020sparse}, Algorithm \ref{algorithm1} does not solve the additional optimization problem. Thus it can avoid many unnecessary numerical errors and make the estimates of true zeros exactly zero.

\begin{figure*}[htbp]
	\centering
	\includegraphics[width=0.9\textwidth]{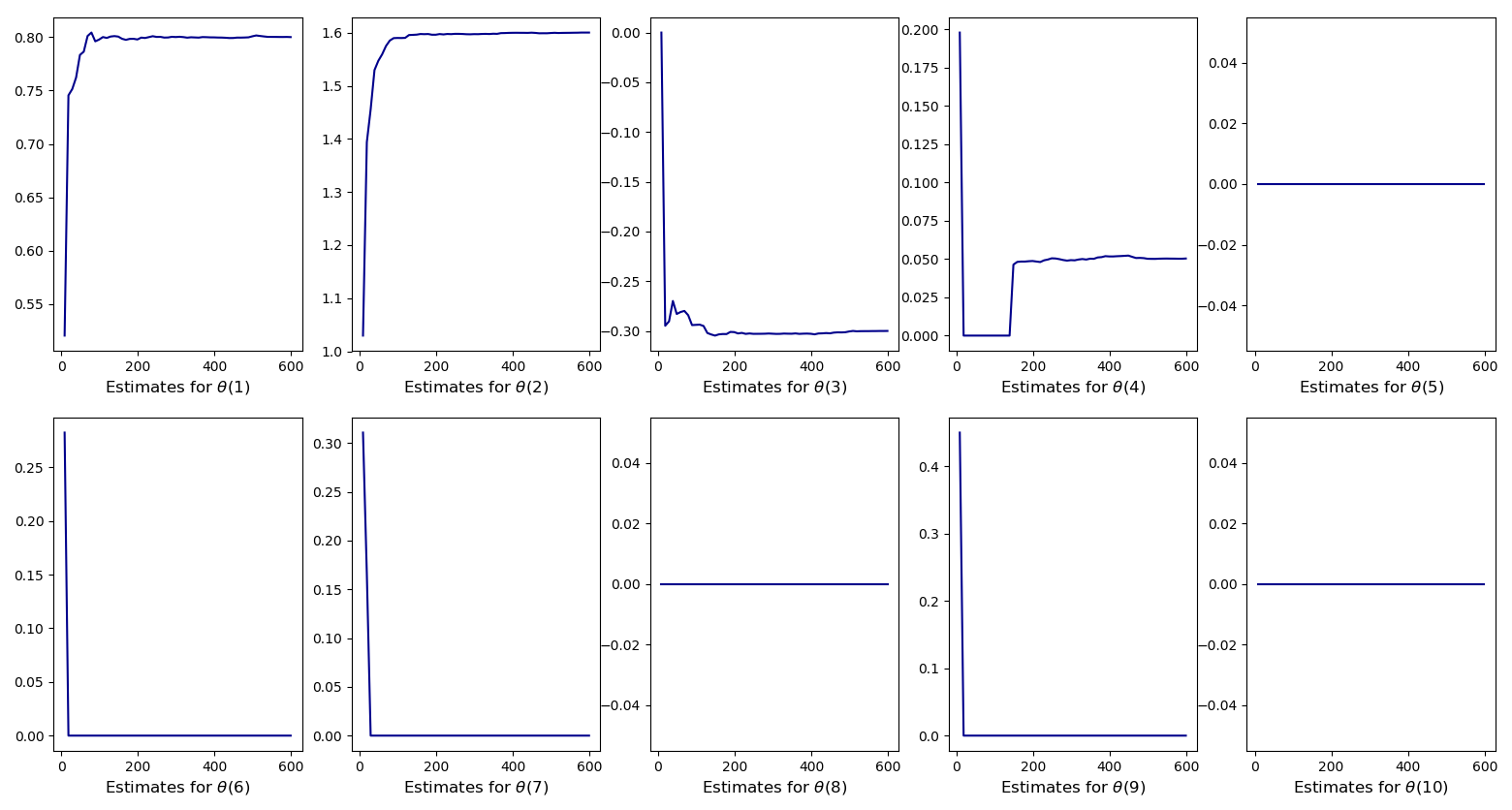}
	\caption{Estimate sequences}
	\label{fig1:sparse}
\end{figure*}

\begin{table*}[htbp]
    \centering
	\caption{Averaged estimates by our method, least squares and LASSO from ten simulations}
	\label{tab1:sparse}
	\tabcolsep=0.5cm
	\resizebox{\textwidth}{!}{
    \begin{tabular}{llllll}
			\toprule
			&N=100 &N=200 &N=300 &N=400 &N=500\\
			\midrule
			Estimates for $\theta(5)$ \\
			By Algorithm \ref{algorithm1}& 0 & 0 & 0 & 0 & 0\\
			By least squares & -3.8403$\times 10^{-3}$ & -1.5291$\times 10^{-3}$ & 1.9044$\times 10^{-4}$ & -7.8872$\times 10^{-5}$ & 1.1098$\times 10^{-4}$\\
			By LASSO & 5.0651$\times 10^{-5}$ & -2.7055$\times 10^{-5}$ & -3.1356$\times 10^{-5}$ & -2.5534$\times 10^{-6}$ & -1.4774$\times 10^{-6}$\\
			By \cite{zhao2020sparse} & 1.8805$\times 10^{-6}$ & -5.3767$\times 10^{-7}$ & -5.8867$\times 10^{-7}$ & -6.3817$\times 10^{-7}$ & -2.4387$\times 10^{-7}$\\
			\midrule
			Estimates for $\theta(6)$& \\
			By Algorithm \ref{algorithm1}& 0 & 0 & 0 & 0 & 0 \\
			By least squares& 4.7193$\times 10^{-3}$ & 1.4462$\times 10^{-3}$ & 8.8818$\times 10^{-4}$ & 8.7922$\times 10^{-4}$ & 7.1250$\times 10^{-4}$\\
			By LASSO & 1.0038$\times 10^{-4}$ & 4.2274$\times 10^{-5}$ & 3.8088$\times 10^{-5}$ & -9.0536$\times 10^{-6}$ & -1.9497$\times 10^{-5}$\\
			By \cite{zhao2020sparse} & 1.9487$\times 10^{-6}$ & 9.5211$\times 10^{-7}$ & -5.2565$\times 10^{-7}$ & -4.8326$\times 10^{-7}$ & -6.0134$\times 10^{-7}$\\
			\midrule
			Estimates for $\theta(7)$& \\
			By Algorithm \ref{algorithm1}& 0 & 0 & 0 & 0 & 0 \\
			By least squares & 1.6445$\times 10^{-3}$ & 1.4848$\times 10^{-3}$ & 1.0246$\times 10^{-3}$ & 8.8689$\times 10^{-5}$ & -3.0212$\times 10^{-5}$\\
			By LASSO & 8.3569$\times 10^{-5}$ & 1.1443$\times 10^{-5}$ & -3.9639$\times 10^{-6}$ & 2.9968$\times 10^{-6}$ & -1.8478$\times 10^{-6}$\\
			By \cite{zhao2020sparse} & -3.8009$\times 10^{-9}$ & 1.4456$\times 10^{-10}$ & -1.3981$\times 10^{-10}$ & 2.8872$\times 10^{-11}$ & -7.8004$\times 10^{-11}$\\
			\midrule
			Estimates for $\theta(8)$& \\
			By Algorithm \ref{algorithm1}& 0 & 0 & 0 & 0 & 0 \\
			By least squares & -3.6471$\times 10^{-3}$ & 3.6763$\times 10^{-3}$ & 1.1463$\times 10^{-3}$ & 9.2123$\times 10^{-4}$ & -9.6761$\times 10^{-5}$\\
			By LASSO & 3.4481$\times 10^{-4}$ & -7.1948$\times 10^{-5}$ & -4.0180$\times 10^{-5}$ & -4.0338$\times 10^{-5}$ & -1.0343$\times 10^{-5}$\\
			By \cite{zhao2020sparse} & -7.8657$\times 10^{-6}$ & -8.4547$\times 10^{-7}$ & -3.2375$\times 10^{-7}$ & -3.5296$\times 10^{-7}$ & -4.4807$\times 10^{-8}$\\
			\midrule
			Estimates for $\theta(9)$& \\
			By Algorithm \ref{algorithm1}& 0 & 0 & 0 & 0 & 0 \\
			By least squares & -4.6465$\times 10^{-3}$ & -2.1687$\times 10^{-3}$ & -8.9557$\times 10^{-4}$ & -4.8013$\times 10^{-4}$ & -6.1303$\times 10^{-5}$\\
			By LASSO & 1.2557$\times 10^{-5}$ & 1.5550$\times 10^{-5}$ & -2.0270$\times 10^{-5}$ & -2.7999$\times 10^{-5}$ & -3.2130$\times 10^{-6}$\\
			By \cite{zhao2020sparse} & 1.1125$\times 10^{-8}$ & 1.4665$\times 10^{-8}$ & -3.8394$\times 10^{-9}$ & -3.5668$\times 10^{-9}$ & -1.4756$\times 10^{-9}$\\
			\midrule
			Estimates for $\theta(10)$& \\
			By Algorithm \ref{algorithm1}& 0 & 0 & 0 & 0 & 0 \\
			By least squares & 3.5405$\times 10^{-3}$ & 1.0033$\times 10^{-3}$ & 4.1614$\times 10^{-4}$ & -4.5427$\times 10^{-4}$ & -1.3878$\times 10^{-4}$\\
			By LASSO & 7.0112$\times 10^{-5}$ & 1.6159$\times 10^{-4}$ & 7.9057$\times 10^{-5}$ & 5.5963$\times 10^{-5}$ & -1.0325$\times 10^{-5}$\\
			By \cite{zhao2020sparse} & -1.9684$\times 10^{-5}$ & 3.7883$\times 10^{-5}$ & -3.8752$\times 10^{-6}$ & -2.1050$\times 10^{-6}$ & -7.3450$\times 10^{-6}$\\
			\bottomrule
	\end{tabular}}
\end{table*}
\end{example}

\section{Conclusion}\label{conclu}
In this paper, the identification of sparse unknown parameters in stochastic regression models is studied. We proposed an efficient algorithm for generating sparse estimates without any penalty terms. Under the weakest non-persistent excited condition, we proved that the set of zero elements in the unknown sparse parameter vector could be correctly identified in finite observations, and the non-zero elements almost surely converge to the true values. For the Hammerstein system, we give the determined number of finite observations. Some meaningful problems deserve further consideration,
e.g., the application of the proposed algorithm to  identify time-varying unknown sparse parameters.
The distributed algorithm to estimate the unknown parameter using local measurement, the adaptive control by  using the estimation algorithm based on the sampled data, and the continuing execution of the algorithm when the set of zero elements is identified in finite observations, i.e., the influence of the initial value on the algorithm.

                                                   







\end{document}